\documentclass{amsart}
\usepackage{graphicx} 
\usepackage{amsmath,amsfonts,amssymb,enumerate,url}
\usepackage{amsthm}
\usepackage{tikz}
\usepackage[sort,nocompress]{cite}
\usetikzlibrary{patterns}
\usetikzlibrary{decorations}
\usetikzlibrary{decorations.markings}
\usepackage{comment}
\newtheorem{theorem}{Theorem}[section]
\newtheorem{lemma}[theorem]{Lemma}
\newtheorem{corollary}[theorem]{Corollary}

\newtheorem{proposition}[theorem]{Proposition}
\theoremstyle{definition}
\newtheorem{definition}[theorem]{Definition}
\newtheorem{notation}[theorem]{Notation}

\theoremstyle{remark}

\def\final{0}  
\def\iflong{\iffalse}
\ifnum\final=0 
\newcommand{\Tamas}[1]{{\color{red}[{\tiny \textbf{Tamas:} \bf #1}]\marginpar{\color{red}*}}}
\newcommand{\student}[1]{{\color{orange}[{\tiny \textbf{Student:} \bf #1}]\marginpar{\color{orange}*}}}
\else 
\newcommand{\Tamas}[1]{}
\newcommand{\student}[1]{}
\fi
\newcommand{\punct}[1]{\langle #1\rangle}

\newcommand{\R}{\ensuremath{{\mathbb{R}}}}

\newcommand{\N}{\ensuremath{{\mathbb{N}}}}

\newcommand{\col}

\newcommand{\Iso}{\ensuremath{{\mathrm{Isom}^+(\R^2)}}}
\newcommand{\IsoS}{\ensuremath{{\mathrm{Isom}^+(S)}}}
\newcommand{\eps}{\varepsilon}
\newcommand{\sm}{\setminus}
\newcommand{\su}{\subset}
\newcommand{\de}{\delta}
\newcommand{\Pt}{\tilde{P}}
\newcommand{\Qt}{\tilde{Q}}

\newcommand{\Et}{\tilde{E}}
\DeclareMathOperator*{\concat}{Concat}

\DeclareMathOperator*{\inte}{int}

\DeclareMathOperator{\id}{identity}

\begin{document}

\title[Continuous Kakeya type problems with special points]{Continuous Kakeya type problems with special points that sweep much more area than the rest} 
\author[B.~Joseph, T.~Keleti, M.~K\"ok\'enyesi, A.~M\'ath\'e and A.~Schmidt]{Benjamin Joseph, Tam\'as Keleti, M\'ark K\"ok\'enyesi, Andr\'as M\'ath\'e and Andrew Schmidt}

\address
{Department of Mathematics and Statistics,
Warner Hall,
Middlebury College,
Middlebury, VT 05753, USA}

\email{bjoseph@middlebury.edu}

\address
{Institute of Mathematics, ELTE E\"otv\"os Lor\'and University,
H-1117, 
Budapest, Hungary
and \newline\indent
HUN-REN Alfr\'ed R\'enyi Institute of Mathematics,
H-1053, Budapest, Hungary}

\email{tamas.keleti@gmail.com}

\address
{Department of Pure Mathematics and Mathematical Statistics, University
of Cambridge, Cambridge, CB3 0WB, UK}

\email{mark.p.kokenyesi@gmail.com}

\address
{Mathematics Institute, University of Warwick, Coventry, CV4 7AL, UK}

\email{a.mathe@warwick.ac.uk}

\address
{Department of Mathematics, Hillsdale College, Hillsdale, Michigan, USA}

\email{ascm2023@gmail.com}

\thanks{The authors would like to thank the Budapest Semesters in Mathematics program
for providing the framework under which part of this research was conducted.
The second author was also supported  
by the Hungarian National Research, Development and Innovation Office -- NKFIH, 146922.
This material is partly based upon work supported by the Swedish Research Council under grant no. 2021-06594 while the second author was in residence at Institut Mittag-Leffler in Djursholm, Sweden during the Fall semester of 2026.
}

\begin{abstract}
For any $\eps>0$ we construct a continuous motion of a line segment in the plane, such that the center point of the line segment traverses a given bounded set except for a subset of area less than $\eps$, while the rest of the line segment sweeps only an area less than $\eps$. 
We further show that a circle with a pair of antipodal points removed has the \emph{Strong Kakeya Property} in the plane: that is, given any two copies of such a punctured circle, we can continuously move one to the other, while staying in a set of arbitrarily small Lebesgue measure.
We give a full characterization of those subsets of a great circle that have this property on the sphere.
We also show how the first mentioned result can be applied to obtain discontinuity results about the area swept by continuously moving objects.
\end{abstract}

\maketitle

\section{Introduction}

\subsection{Historical background and the main results}
To answer Kakeya's needle problem,  Besicovitch \cite{Be} proved in 1928 that a unit line segment can be fully rotated continuously in the plane in arbitrarily small area.
Following \cite{CHL}, we say that a set $E\subset \R^2$ has the \emph{Strong Kakeya Property} if for any two isometric copies, we can move from one to the other by a continuous rigid motion, so that the set sweeps out arbitrarily small area (see the precise definition in Section~\ref{s:prel}). 
It is easy to show (using the so-called P\'al joins) that the above result of Besicovitch implies that 
the line segment has the strong Kakeya property. In \cite{CHL} it is shown that any closed and connected set with the strong Kakeya property is contained in a line or a circle. 

Cunningham \cite{Cu} asked in 1974 whether circular arcs have the strong Kakeya property - without using this terminology. It was proved by H\'era and Laczkovich \cite{HL} that the answer is affirmative, at least for circular arcs of small enough angle.
This was significantly extended by Chang and Cs\"ornyei \cite{CC}, who showed that a circle with a pair of open antipodal arcs removed also has the strong Kakeya property. 
One of our main results extends this even further: we prove that in fact, it is sufficient to remove only two antipodal points.

\begin{theorem}\label{t:antipodal-intro}
The circle with two antipodal points removed has the strong Kakeya property.
\end{theorem} 



In \cite{Cu} Cunningham also studied other closely related problems. 
Cunningham defines a \emph{bird} as a line segment, which consists of a small central segment, called the \emph{body}, and the remaining two long end segments, called the \emph{wings}. 
He proved that a bird can be moved continuously in the plane in such a way that the body touches every point of a given bounded set, while the wings sweep over a set of arbitrarily small area.
We also extend this result further by showing that almost the same is true even if the body of the bird is a single point: in this case the body can touch an arbitrarily large part of the given bounded set, while the wings sweep a set of arbitrarily small area.

\begin{theorem}\label{t:bird-intro}
    Let $E$ be a line segment, $q$ be any point of $E$,
$S$ be any bounded subset of the plane and $\eps>0$. 
The line segment $E$
can be moved continuously so that during the motion $E\sm\{q\}$ sweeps an area less than $\eps$, while the single point $q$
sweeps all points of $S$ but a set of area less than $\eps$.
\end{theorem}

In \cite{Cu} Cunningham also solved the spherical version of the Kakeya needle problem: he found the greatest lower bound $k_a$ of surface areas of sets on the unit sphere in which a closed great circle arc of angle $a$ can be continuously turned around. He proved that $k_a=0$ for $a<\pi$, $k_a=2\pi$ for $\pi\le a<2\pi$ and of course $k_{2\pi}=4\pi$.

The notion of strong Kakeya property naturally extends to the unit sphere.
The above $k_a=0$ result for $a<\pi$ implies that closed great circle arcs of angle less than $\pi$ have the strong Kakeya property. 
We extend this result to the great circle with two antipodal points removed.
Remarkably, combining this result with a simple argument in \cite{Cu} leads to a complete characterization of those subsets of a great circle that have the strong Kakeya property on the sphere.

\begin{corollary}\label{c:sphere-char-intro}
A subset $A$ of a great circle $E$ has the strong Kakeya property on the sphere if and only if $E\setminus A$ contains a pair of antipodal points.
\end{corollary}

In the plane a halfline cannot have the strong Kakeya property since we need to cover infinite area to move it to a different line. 
But this trivial obstacle disappears if instead of Lebesgue measure we consider a fixed absolutely continuous finite Borel measure $\mu$.
Clearly, unless $\mu=0$, a full line still does not have the strong Kakeya property with respect to $\mu$. 
But we could show that any proper subset of a line does:

\begin{theorem}\label{t:punched-line-intro}
A line punctured at one point has the strong Kakeya property in the plane with respect to any absolutely continuous finite Borel measure.
\end{theorem} 

We also show applications of Theorem~\ref{t:bird-intro} to obtain some discontinuity results about the area swept by continuously moving objects.
Using a construction of Talagrand \cite{Ta}, Fleiner \cite{Fl} proved that the plane can be moved continuously such that the path of the origin has positive Lebesgue measure but the path of any other point in $[-1,1]^2$ has measure zero,
thus the area swept by a point does not depend continuously on the point.
Theorem~\ref{t:bird-intro} will easily give another construction for this noncontinuity, and with a bit more effort it will also give the following noncontinuity result about the area swept by line segments.

\begin{theorem}\label{t:opendiagonal-intro}
 If $E\su\R^2$ is an open line segment and $v\in\R^2$ is 
 not perpendicular
 to $E$, then the plane can be moved continuously so that the area swept by $E+tv$ is not continuous in $t$. 
\end{theorem}

\subsection{About the proofs and the structure of the paper}
All of the above mentioned main results follow from the same construction.
This construction uses heavily a theorem of Chang and Cs\"ornyei \cite{CC}.
In Section~\ref{s:prel}, after fixing the terminology, we prove some useful lemmas, state the above mentioned theorem of Chang and Cs\"ornyei, and prove its direct corollaries. 
Section~\ref{s:construction} is the most technical part of the paper, it contains the constructions obtained from the corollaries of Section~\ref{s:prel}.
Then in Section~\ref{s:main} we state precisely and prove Theorems~\ref{t:antipodal-intro}, \ref{t:bird-intro} and \ref{t:punched-line-intro} using the constructions of Section~\ref{s:construction}, and we obtain Corollary~\ref{c:sphere-char-intro} from Theorem~\ref{t:punched-line-intro}.
Finally, Section~\ref{s:applications} contains the applications about the discontinuity of the area swept by continuously moving objects, including the precise statement and the proof of Theorem~\ref{t:opendiagonal-intro}.

\subsection{AI declaration}
All the results and proofs were found and written up by the authors.
ChatGPT-5.6 and ChatGPT-6 helped to find errors in earlier versions 
and to find useful references.

\section{Preliminaries}\label{s:prel}

For $X\subset \R^k,$ we let $|X|$ denote the $k$-dimensional Lebesgue (outer) measure of $X$.  
When we work on a sphere then $|X|$ denotes the surface outer measure of $X$.

Let $\Iso$ denote the space of orientation-preserving isometries of $\R^2$.
Every $\rho\in\operatorname{Isom}^+(\mathbb R^2)$ can be written
uniquely as
\[
\rho(x)=R_\theta x+v,
\qquad \theta\in\mathbb R/2\pi\mathbb Z,\quad v\in\mathbb R^2.
\]
We equip $\operatorname{Isom}^+(\mathbb R^2)$ with the complete metric
\[
d_{\mathrm{Iso}}(\rho,\rho')
=d_{\mathbb S^1}(\theta,\theta')+|v-v'|,
\]
where $d_{\mathbb S^1}$ denotes angular distance modulo $2\pi$.
We will work with continuous paths $P:[a,b]\to\mathrm{Isom}^+(\R^2).$ 
For paths with a common domain let $d$ be the supremum metric with respect to the fixed metric on \Iso.
This is clearly also a complete metric.  

We define a \emph{translational segment} as a path of isometries 
$$
P:[a,b]\to\mathrm{Isom}^+(\R^2), \qquad 
P(t):x \mapsto \rho(x)+(t-a)\vec{v}
$$ 
for some fixed $\rho\in\Iso$ and $\vec{v}\in\R^2$. 
We define a \emph{rotational segment} as a path of isometries 
$P:[a,b]\to\mathrm{Isom}^+(\R^2)$ such that for some fixed 
$\rho\in\Iso$, $z\in \R^2$ and $s\in \R $, $P(t)\circ \rho^{-1}$ is a rotation of $(t-a)s$ radians around $z$. 
We define a \emph{segment of motion} as either of the above types of paths,
and a \emph{polygonal path} as a (continuous) concatenation of a finite sequence of segments of motion, for which we will use the notation $P=\concat_{i=1}^n L_i$. 
Note that with the above defined metric $d_{\mathrm{Iso}}$, every polygonal path is Lipschitz.
By a \emph{path} we always mean a continuous
path $P:[a,b]\to \Iso$  and we call it a \emph{path from $\rho$ to $\rho'$} (or \emph{with endpoints $\rho$ and $\rho'$}) if $P(a)=\rho$ and $P(b)=\rho'$.

For a path $P:[a,b]\to \Iso$, and a set $E\subset\R^2,$ we will denote the set swept by $E$ under $P$ as $P[E]$; that is, $$P[E]=\bigcup_{t\in[a,b]}P(t)(E).$$

A set $E\su\R^2$ is said to have the \emph{Strong Kakeya Property} if for any $\rho_1,\rho_2\in\Iso$ and $\eps>0$ there exists a path $P:[a,b]\to\Iso$ such that 
$P(a)=\rho_1$, $P(b)=\rho_2$ 
and $|P[E]|<\eps$.

\begin{lemma}\label{l:Federer}
For some $k\in\N^+$ let $X\su\R^k$ be compact such that its boundary has Lebesgue measure zero and let $C_k(X)$ denote the set of continuous functions $g:X\to\R^k$ with the supremum metric $d_\infty$.
Then for any Lipschitz function $f:X\to\R^k$ and $\eps>0$ there exists an $n\in\N^+$ such that whenever $g\in C_k(X)$ with $d_\infty(f,g)<1/n$ then 
$|f(X)\sm g(X)|<\eps$.
\end{lemma}

\begin{proof}
A $y\in f(X)$ is called a \emph{stable value of $f$} if $y$ is in the range of every $g\in C_k(X)$ sufficiently close to $f$.
By a classical result of Federer, stated explicitly in
\cite[p.~319]{Fe52} and proved, in an equivalent formulation, in the
proof of \cite[Theorem~6.12, pp.~461--462]{Fe46}, almost every
$y\in f(X)$ is a stable value of $f$.

Let $Y$ be the set of stable values of $f$ and let
$$
Y_n=\{y\in Y \colon g\in C_k(X), d_\infty(g,f)<1/n \Rightarrow y\in g(X)\}.
$$
Each $Y_n$ is closed and by definition, $Y_n\su Y_{n+1}$ and $Y=\bigcup_{n} Y_n$.
By the above mentioned result of Federer, $|f(X)\sm Y|=0$, so there exists an $n$ such that $|f(X)\sm Y_n|<\eps$.
Since by definition $Y_n\su g(X)$ for any $g\in C_k(X)$ with $d_\infty(f,g)<1/n$, this completes the proof.
\end{proof}

The first part of the following 
lemma was a key tool in \cite{CC}.

\begin{lemma}\label{uppersemi}
Let $P$ be a path, 
$E \su \R^2$ be compact and $\eps>0$.
\begin{enumerate}[(a)]
 \item \label{part:upper}
 There exists some $\delta>0$ such that for any path $Q$ within distance $\delta$ of $P$, we have 
    \begin{align*}
    |Q[E]\sm P[E]| < \eps.
 \end{align*}
 \item \label{part:lower}
 If $P$ is a Lipschitz path and $E$ admits a Lipschitz parametrization $\gamma:[0,1]\to E$ onto $E$ then there exists some $\delta>0$ such that for any path $Q$ within distance $\delta$ of $P$, we have 
    \begin{align*}
    |P[E]\sm Q[E]| < \eps.
 \end{align*}
In particular, this holds if $P$ is a polygonal path and $E$ is a closed line segment or a circle.

\end{enumerate}     
   
\end{lemma}

\begin{proof}
Part (\ref{part:upper}) is \cite[Lemma 2.1]{CC}.
To prove (\ref{part:lower})
let $[a,b]$ be the domain of $P$ and $Q$ and $X=[a,b]\times [0,1]$.
For $(u,v)\in X$ let $f(u,v)=P(u)(\gamma(v))$ and $g(u,v)=Q(u)(\gamma(v))$.
Put $M=\max\{1,\sup_{x\in E}|x|\}$. From the definition of
$d_{\mathrm{Iso}}$ for any $\rho, \rho'\in\Iso$ and $x\in E$ we have
$|\rho(x)-\rho'(x)|\le M d_{\mathrm{Iso}}(\rho,\rho')$.
Hence $f$ is Lipschitz and $d_\infty(f,g)\le M d(P,Q)$, so
Lemma~\ref{l:Federer} gives the claim. 
\end{proof}

We use the following classical result of Cunningham \cite{Cu}:
 
\begin{theorem}[Cunningham Bird Theorem]\label{Cunningham}
    Let $E$ be a line segment, and $U_0$ an open subsegment of $E$. For any bounded set $S\su\R^2$ and $\varepsilon>0,$ there is a polygonal path $P$ so that $S\su P[U_0]$, but $|P[E\setminus U_0]|<\varepsilon.$ 
\end{theorem}

(Although Cunningham does not state explicitly that this path is polygonal, this follows from the construction in his proof: the motion is obtained as a finite concatenation of translations and rotations. Thus the path in Theorem~\ref{Cunningham} can be taken to be polygonal.) We use the following theorem from \cite{CC} (Theorem 5.1 and Remark 5.11) for all of our constructions.

\begin{theorem}[Chang-Cs\"ornyei]\label{CClemma}
    Let $E$ be a smooth curve of finite length, $\varepsilon>0$ and $\rho\in \mathrm{Isom}^+(\mathbb{R}^{2})$ be arbitrary. If $\rho$ is a rotation, we take any line $\ell$ through the center of the rotation; 
if $\rho$ is a translation, we take any line $\ell$ perpendicular to the translation. Then there exists a sequence of subsets of $\ell$, $W_{1},W_{2},\cdots,W_{n}$ where either $W_{i}$ is an $\eps$-segment of $\ell$, or is $\ell \setminus B_{1/\varepsilon}((0,0))$ with the following properties. Let $\nu_{x}$ be the line normal to $E$ at $x.$ 
For each $W_{i}$ which is a line segment, we define 
$$Y_{i}=\{ x\in E: \nu_{x}\cap W_{i}\ne \emptyset\},$$
and for the other $W_{i}$ we define 
$$ Y_{i}=\{ x\in E: \nu_{x}\cap W_{i}\ne \emptyset \}\cup \{ x\in E:\nu_{x} \text{ parallel to } \ell \}.$$
Then in any neighborhood of any prescribed segment of motion from identity to $\rho$, there exists a polygonal path $P=\concat_{i=1}^nL_{i},$ from the identity to $\rho$, where $L_{i}$ are segments of motion, so that $$ \left\vert\bigcup_{i=1}^n L_{i}[E\setminus Y_{i}] \right\vert<\varepsilon.$$
\end{theorem}

For our purposes, $E$ will always be a line segment or a circle. In these subcases, we get the following corollaries. 

\begin{corollary}\label{cor 5}
Let $E$ be a line segment, $\rho_1, \rho_2\in \mathrm{Isom}^+(\mathbb{R}^{2})$ and $\varepsilon>0$ be arbitrary. 
Then in any neighborhood of any prescribed segment of motion from $\rho_1$ to $\rho_2$, there exists a polygonal path $P=\concat_{i=1}^n L_{i}$ from $\rho_1$ to $\rho_2$ and points $u_{i}\in E$ so that $$ \left|\bigcup_{i}L_i[E\setminus B_{\varepsilon}(u_{i})]\right|<\varepsilon.$$
\end{corollary}
\begin{proof}
   Without loss of generality we can suppose that $\rho_1$ is the identity. Let $\rho=\rho_2$. We want to show that with the right choice of $\varepsilon$ and $\ell$ in Theorem \ref{CClemma}, we can make the $Y_{i}$ equal to an $\eps$-segment of our line segment $E$. 
If $\rho$ is a translation parallel to $E$ then the statement of the corollary is trivial. 
Otherwise we can choose $\ell$ to be not perpendicular to $E$.
Then
$Y_{i}$ is the intersection of $E$ and the projection of $W_{i}$ to the line containing $E$. We can choose $\varepsilon$ small enough so $\ell \setminus B_{1/\varepsilon}$ does not contain any points that project onto the segment $E$. 
%
\end{proof}

\begin{corollary}\label{cor 6}
    Let $E$ be a circle, $\rho_1, \rho_2\in \mathrm{Isom}^+(\mathbb{R}^{2})$ and $\varepsilon>0$ be arbitrary. 
Then in any neighborhood of any prescribed segment of motion from $\rho_1$ to $\rho_2$, there exists a polygonal path $P=\concat_{i=1}^n L_{i}$ from $\rho_1$ to $\rho_2$,    
and antipodal pairs of points $z_{i},z_i'\in E$ so that $$ \left|\bigcup_{i}L_i[E\setminus(B_{\varepsilon}(z_{i})\cup B_\eps(z_{i}'))]\right|<\varepsilon.$$ 
\end{corollary}

\begin{proof}
Again, we can suppose that $\rho_1$ is the identity and we want to apply Theorem~\ref{CClemma} for $\rho=\rho_2$.
    If $\rho$ is a rotation around the center of our circle, the claim follows immediately with $P$ a rotational segment. Otherwise, choose $\ell$ in Theorem $\ref{CClemma}$ so that it does not contain the center of the circle.

    If $W_i=\ell\setminus B_{1/\varepsilon}((0,0))$, then $Y_i$ is a pair of antipodal arcs. By choosing $\varepsilon$ small enough, we can make these arcs arbitrarily small. Since $\ell$ does not contain the center of the circle, if $W_i$ is a segment, then $Y_i$ is also a pair of antipodal arcs, which can be made arbitrarily small by choosing $\eps$ small enough.
\end{proof}

\section{Construction}\label{s:construction}

Chang and Csörnyei used Theorem~\ref{CClemma} to take a limit and prove that by deleting a single point at every moment the remaining points sweep measure zero. We rebuild a similar argument, showing we can sacrifice the measure zero but gain continuity of the point being deleted while still keeping the measure small.

\begin{notation}
If $E$ is a circle or a line segment then for \(p\in E\), let
\[
\punct{p}=
\begin{cases}
\{p\}, & \text{if \(E\) is a line segment},\\
\{p,p^*\}, & \text{if \(E\) is a circle},
\end{cases}
\]
where \(p^*\) denotes the point antipodal to \(p\).
\end{notation}

We begin with a technical statement, which says that we can iterate the corollaries of Theorem~\ref{CClemma}, so that the 
removed sets
nest. 

\begin{lemma}\label{sequence}
Let $E$ be a closed line segment (or circle), $U_0$ an open interval of $E$ (or a pair of antipodal open arcs), sequences of positive reals $(\varepsilon_n),(\eta_n),(\delta_n)$, and a polygonal path $P_0:[0,1]\to\Iso$ be given. 
Then there exist polygonal paths $P_n:[0,1]\to\Iso$ comprised of segments of motion $L_v:I_v\to\Iso$ for $v\in V_n$, where $V_n$ is a finite set of integer vectors of length $n+1$, along with nonempty open subsets $U_v$ of $E$, for $v\in V_n$ and $\mu_n\in(0,\eta_n)$
with the following properties.
\begin{enumerate}
	\item \label{c:endpoints}
	For all $n>0$ and $v \in V_{n-1}$, the restriction $P_n\restriction I_v$ is a polygonal path with the same endpoints as the segment of motion $L_v$.
    \item \label{c:deltan} 
    For all $n>0$ and $v \in V_{n}$, there is some $p\in E$ so that $U_v\subset B_{\delta_n}(\langle p\rangle)$ (that is, $U_v$ is covered by a segment (or pair of antipodal arcs) each of diameter at most $2\delta_n$).
    \item \label{c:chopv}
    For all $v\in V_0$ we have $U_v=U_0$, and for all $n>0$ and $v \in V_n$ we have $v^-\in V_{n-1}$, $U_v \su U_{v^-}$ and $I_v \su I_{v^-}$, where $v^-$ denotes the vector obtained by deleting the last element of $v$.
    \item \label{c:dPmPn}
    For all $n>m \geq 0$,
    \begin{align*}
        d(P_m,P_n) < \mu_m<\eta_m.
    \end{align*}
    \item \label{c:P0EminusPEsmall}
For all $n\ge 1$,    
    \begin{align*}
        |P_0[E]\setminus P_n[E]| < \sum_{i=1}^n \varepsilon_i.
    \end{align*}
    \item \label{c:unionLvEminusUvsmall}
For all $n\ge 1$, 
    \begin{align*}
        \left|\bigcup_{v \in V_n}L_v[E\backslash U_v]\right| < \left|P_0[E\setminus U_0]\right| + \sum_{i=1}^{n}\varepsilon_i.
    \end{align*}
    \item \label{c:Q}
    For all $n\ge 1$ and any path $Q:[0,1]\to\Iso$ such that $d(Q,P_n)\le\mu_n$,
    \begin{align*}
        \left|\bigcup_{v \in V_n}(Q\restriction I_v)[E\backslash U_v]\right| < \left|P_0[E\setminus U_0]\right| + \sum_{i=1}^{n}\varepsilon_i.
    \end{align*}

\end{enumerate}
\end{lemma}
\begin{proof}

The proof is by induction. Let $V_0$ index the segments of motion $L_v$ which compose $P_0$, and for each $v \in V_0$, let $U_v = U_0$. 
Let $\mu_0\in (0,\eta_0)$.
Then conditions \ref{c:endpoints}-\ref{c:Q} clearly hold for $n=0$, so our induction hypothesis is that the result is true for a finite sequence of $n$ paths, $P_0,\ldots,P_{n-1}$, with $n \geq 1$. 

We construct $P_n$ by first constructing polygonal paths $P_n^v:I_v\to\Iso$ corresponding to each $v \in V_{n-1}$, where $P_n^v$ comes with a corresponding index set $V_n^v$. 
For each $v\in V_{n-1}$, by applying Corollary~\ref{cor 5} or Corollary~\ref{cor 6} to $L_v$, we obtain
a finite index set $V_n^v\su \{v\}\times \N$, polygonal paths
\begin{align*}
    P_n^v = \concat_{w \in V_n^v} L_w
\end{align*}
with the same endpoints as $L_v$,
and for each $w\in V_n^v$ a point $x_w\in E$ such that 
\begin{align}\label{e:unionLw}
    \left|\bigcup_{w \in V_n^v} L_w[E\backslash B_{\delta_n}(\langle x_w\rangle ) ]\right| < \frac{\varepsilon_n}{2|V_{n-1}|},
\end{align} 
where for index sets $|\cdot|$ denotes cardinality.

Let $V_n=\cup_{v\in V_{n-1}} V_n^v$ and $P_n=\concat_{v\in V_{n-1}}P_n^v$.
By construction, condition \ref{c:endpoints} holds for $n$.
Since Corollaries \ref{cor 5} and \ref{cor 6} guarantee each $P_n^v$ can  be constructed to be arbitrarily close to $L_v$, we can make $P_n$ arbitrarily close to $P_{n-1}$. 
Since, by the induction hypothesis, $d(P_m,P_{n-1})<\mu_m$ for all $m < n-1$,
we can choose $0 < \mu_{n-1}' < \mu_{n-1}$ small enough that any path within distance $\mu_{n-1}'$ of $P_{n-1}$ is also strictly within distance $\mu_m$ of $P_m$ for all $m < n-1$. 
We construct $P_n$ to be within distance $\mu_{n-1}'$ of $P_{n-1}$.
Thus condition \ref{c:dPmPn} holds for $n$ (and every $m<n$).

 Additionally, by part \eqref{part:lower} of  Lemma~\ref{uppersemi}, constructing $P_n$ close enough to $P_{n-1}$ we can guarantee
\begin{equation*}
|P_{n-1}[E]\sm P_n[E]|<\eps_n. 
\end{equation*}  
This and the induction hypothesis give condition \ref{c:P0EminusPEsmall}.

Similarly, for each $v\in V_{n-1}$, constructing $P_n^v$ close enough to $L_v$ we can also guarantee by part \eqref{part:upper} of  Lemma~\ref{uppersemi},
\begin{equation*}
|P_n^v[E\sm U_v]\sm L_v[E\sm U_v]|< \frac{\eps_n}{2|V_{n-1}|}.
\end{equation*}
Using this and the induction hypothesis we obtain
\begin{equation}\label{e:88}
\begin{aligned}
    \left|\bigcup_{w \in V_n}L_w[E \backslash U_{w^-}]\right| =
    \left|\bigcup_{v \in V_{n-1}}P_n^v[E\backslash U_{v}] \right|
    &< \left|\bigcup_{v \in V_{n-1}}L_v[E\backslash U_{v}] \right| + \frac{\varepsilon_n}{2} \\
    &\le |P_0[E\setminus U_0]| + \sum_{i=1}^{n-1}\varepsilon_i + \frac{\varepsilon_n}{2}.
\end{aligned}
\end{equation}

We define $U_w = B_{\delta_n}(\langle x_w\rangle)\cap U_{w^-}$, if the latter is nonempty. 
Otherwise, 
define $U_w=B_{\delta_n}(\langle x\rangle)\cap U_{w^-}$, where $x$ is chosen so that $U_w$ is nonempty.
Thus, conditions~\eqref{c:deltan} and \eqref{c:chopv} also hold for $n$, so it 
remains to guarantee conditions~\eqref{c:unionLvEminusUvsmall} and \eqref{c:Q}.

Observe that using \eqref{e:unionLw} we obtain
\begin{equation}\label{e:89}
\begin{aligned}
    \left|\bigcup_{w \in V_n}L_w[E\backslash B_{\delta_n}(\langle x_w\rangle ) ]\right| &\leq \sum_{v \in V_{n-1}} \left|\bigcup_{w \in V_n^v}L_w[E\backslash B_{\delta_n}(\langle x_w\rangle)]\right| \\
    &< |V_{n-1}|\cdot \frac{\varepsilon_n}{2|V_{n-1}|} 
    = \frac{\varepsilon_n}{2}.
\end{aligned}
\end{equation}
Using that for any $w\in V_n$, 
$B_{\delta_n}(\langle x_w\rangle)\cap U_{w^-}\su U_w$, and also \eqref{e:88} and \eqref{e:89} we get
\begin{align*}
    \left|\bigcup_{w \in V_n}L_w[E \backslash U_w]\right| &\leq \left|\left(\bigcup_{w \in V_n}L_w[E\backslash B_{\delta_n}(\langle x_w\rangle )]\right)\cup \left(\bigcup_{w \in V_n}L_w[E \backslash U_{w^-}]\right)\right| \\
    &\leq \left| \bigcup_{w \in V_n}L_w[E\backslash B_{\delta_n}(\langle x_w\rangle )]\right| + \left|\bigcup_{w \in V_n}L_w[E\backslash U_{w^-}]\right| \\
    &< \frac{\varepsilon_n}{2} + \left(|P_0[E\setminus U_0]| + \sum_{i=1}^{n-1}\varepsilon_i +\frac{\varepsilon_n}{2} \right) 
    = \left|P_0[E\setminus U_0]\right| + \sum_{i=1}^{n}\varepsilon_i,
\end{align*}
which gives condition~\ref{c:unionLvEminusUvsmall}.

Recall that $L_v=P_n\restriction I_v$ ($v\in V_n$).
By part \eqref{part:upper} of Lemma~\ref{uppersemi},
by choosing $\mu_n\in(0,\eta_n)$ small enough,
we can make each $|(Q\restriction I_v)[E\sm U_v]\sm L_v[E\sm U_v]|$
as small as we wish simultaneously for all paths $Q:[0,1]\to\Iso$ such that $d(Q,P_n)\le \mu_n$.
Thus, using condition~\ref{c:unionLvEminusUvsmall} we can guarantee condition~\ref{c:Q} as well,
which completes the proof.
\end{proof}

Now, we can consider the limit behavior of such a sequence. Using the properties established above, we can prove that the limit of such a sequence exists, and 
after reparametrization we can make the single deleted point in the limit continuous.

\begin{lemma}\label{lemma:P exists}
    Let $E$ be a closed line segment (or circle), $U_0$ an open interval of $E$ (or pair of antipodal arcs). 
Given some fixed $\varepsilon > 0$, $\eta>0$ and polygonal path $P_0:[0,1]\to\Iso$, there exists a path $\Pt:[0,1]\to\Iso$ with the same endpoints as $P_0$
within distance $\eta$ from a reparametrization of $P_0$,
and a continuous map $f:[0,1]\to \overline{U_0}$ 
(in the circular case $f:[0,1]\to E$) 
such that the following hold:
\begin{equation}\label{e:punchedwithf}
\left|\bigcup_{t \in [0,1]}\Pt(t)(E\backslash\punct{f(t)})\right| < |P_0[E\setminus U_0]| +\varepsilon
\end{equation}
and
\begin{equation}\label{e:P0EminusPE}
\left|P_0[E]\setminus \Pt[E]\right| <  \varepsilon.
\end{equation}

\end{lemma}

\begin{proof}
We apply Lemma \ref{sequence} to obtain a sequence of index sets $V_n$, for every $v\in V_n$ nonempty open sets $U_v\su E$ and segments of motion $L_v:I_v\to \Iso$ and a sequence of polygonal paths $P_n=\concat_{v\in V_n} L_v$ with the properties described in Lemma~\ref{sequence}. 
By reparameterizing each $L_v$ and $P_n$ simultaneously by a single increasing homeomorphism of $[0,1]$ we can guarantee that all the endpoints of all the intervals $I_v$ are dyadic rational numbers.
By statement~\ref{c:dPmPn} of Lemma~\ref{sequence},
if we choose $\eta_n$ converging to zero, then the sequence $(P_n)$ is Cauchy (in the supremum metric), so by completeness it converges uniformly to some (continuous) path $P$,
and we also have $d(P_n,P)\le\mu_n<\eta_n$. 
Let $c$ denote the classical Cantor function (devil's staircase) and by reparameterizing simultaneously again, 
let $\Pt=P\circ c$ and for each $n$ let $\Pt_n=P_n\circ c$.
By choosing $\eta_0<\eta$ we have $d(\Pt,\Pt_0)=d(P,P_0)<\eta$, so $\Pt$ is indeed within distance $\eta$ from a reparametrization of $P_0$.
 
By choosing the sequence $(\eps_n)$ such that $\sum_{i=1}^\infty\eps_i <\eps/2$ and $n$ large enough, using $\Pt_n[E]=P_n[E]$, part \eqref{part:upper} of Lemma \ref{uppersemi} and statement \ref{c:P0EminusPEsmall} of Lemma \ref{sequence} we obtain 
\begin{align*}
    \left|P_0[E]\sm \Pt[E]\right| \leq 
    \left|P_0[E]\sm P_n[E]\right| + \left|\Pt_n[E]\sm \Pt[E]\right| < \eps/2 + \eps/2 = \eps,
      \end{align*}
which gives \eqref{e:P0EminusPE}.

    We define our puncture function $f: [0,1] \to E$ as follows. 
Fix $n\ge 0$. Let 
$$
J_v=\overline{c^{-1}(\inte I_v)} \quad (v\in V_n)
\qquad \text{and} \qquad
J^n=\bigcup_{v\in V_n} J_v,
$$
where $\inte$ denotes interior and $\overline{H}$ denotes the closure of the set $H$.
Since by definition $\{I_v \colon v\in V_n\}$ is a decomposition of $[0,1]$ into closed intervals, overlapping only at their endpoints, and the endpoints of every $I_v$ are dyadic rationals, 
$\{J_v \colon v\in V_n\}$ is a finite collection of pairwise disjoint closed subintervals of $[0,1]$ such that $0$ and $1$ are contained in their union $J^n$.
Similarly, for any $v\in V_n$, $\{J_w \colon w\in V_{n+1}, w^-=v\}$ is a finite collection of pairwise disjoint closed subintervals of $J_v$ such that the endpoints of $J_v$ are contained in the union of these subintervals.

Let $\phi_n:J^n \to V_n$ map $t \in J_v$ to the vector $v \in V_n$.
Since the range of $\Pt$ on $J_v$ is the same as the range of $P$ on $I_v$, we have
\begin{equation}\label{e:above}
\bigcup_{t \in J^n}\Pt(t)(E\sm U_{\phi_n(t)})=
\bigcup_{v \in V_n}(\Pt\restriction J_v)[E\sm U_v]=
\bigcup_{v \in V_n}(P\restriction I_v)[E\sm U_v].
\end{equation}
Since $d(P,P_n)\le\mu_n$ and $\sum_i\eps_i<\eps/2$, we can combine \eqref{e:above} with statement \ref{c:Q} of  Lemma~\ref{sequence} to obtain
\begin{equation}\label{e:JnPtn}
\left|\bigcup_{t \in J^n}\Pt(t)(E\sm U_{\phi_n(t)})\right| < 
\left|P_0[E\setminus U_0]\right| + \frac{\eps}2.
\end{equation}


For each $n\ge 1$ and $v\in V_n$, let $x_v$ be such that $U_v\su B_{\delta_n}(\langle x_v\rangle)$. In the line segment case, we will also require $x_v\in U_0$, which is possible since $U_v\su U_0$. In the line segment case
for $t\in J^n$ let $f_{n}(t)$ be $x_{\phi_{n}(t)}$ 
and extend $f_n$ linearly in the intervals of $[0,1]\sm J^n$. We will define $f_n$ in the circular case after a brief observation.

Let $r$ be the radius of $E$ if $E$ is a circle and let $r>0$ be arbitrary if $E$ is a line segment. 
We wish to show that, by choosing $\de_n$ in Lemma~\ref{sequence} so that $\sum\de_n<\infty$ and $\de_n<r/2$ for every $n$, we can guarantee that $f_{n}$ converges uniformly. 
By $\sum\de_n<\infty$, it is enough to show that for any $n\ge 1$ we have 
\begin{equation}\label{e:fnfnplus1}
|f_n-f_{n+1}|\le \de_n+\delta_{n+1}.
\end{equation} 
In the linear case, \eqref{e:fnfnplus1} follows on $J^{n+1}$ by conditions \eqref{c:chopv} and \eqref{c:deltan} of Lemma \ref{sequence}, and the fact that \begin{equation}\label{e:fcenter}
   \emptyset\neq U_{\phi_n(t)}\su B_{\de_n}(\langle f_n(t)\rangle).
\end{equation} Then in the complementary intervals \eqref{e:fnfnplus1} follows from the linear extensions.

We now prove by induction that we can define continuous $f_n:[0,1]\to E$ on the circle so that \eqref{e:fnfnplus1} and \eqref{e:fcenter} hold
and $f_n$ is constant on $J_v$ for every $v\in V_n$. 
For each $v\in V_1$, by \eqref{c:deltan} of Lemma~\ref{sequence}, there is a $p_v\in E$ such that $U_v\subset B_{\de_1}(\langle p_v \rangle)$. Let $f_1=p_v$ on $J_v$ and extend continuously in the complementary intervals. 
Suppose that $n\ge 1$ and the induction hypothesis holds for $f_n$.
Let $t\in J^{n+1}$.
Since \eqref{e:fcenter} holds for $n$,
there is a $p$, depending only on $\phi_{n+1}(t)$, so that 
$$
\emptyset\neq U_{\phi_{n+1}(t)}\su B_{\de_{n+1}}(\langle p \rangle)\cap B_{\delta_n}(\langle f_n(t)\rangle).
$$ 
Then let $f_{n+1}(t)\in \langle p \rangle$ be such that $B_{\de_{n+1}}(f_{n+1}(t))\cap B_{\delta_n}(f_{n}(t))$ is nonempty. 
Then on $J_{\phi_{n+1}(t)}$ \eqref{e:fnfnplus1} and \eqref{e:fcenter} hold
and $f_{n+1}$ is constant.
Thus at the endpoints of each complementary interval of $[0,1]\sm J^{n+1}$ we have \eqref{e:fnfnplus1} and $\de_n+\delta_{n+1}<r$, which clearly imply that $f_{n+1}$
can be extended to each complementary interval as a continuous function 
so that \eqref{e:fnfnplus1} holds everywhere.

Thus in both cases $f_n$ converges uniformly, so it has a continuous limit
$f:[0,1] \to E$, and in the linear case we have $f:[0,1]\to \overline{U_0}.$
It remains to prove \eqref{e:punchedwithf}.

Let $C$ denote the Cantor set. 
By construction, $\Pt$ is constant on each complementary interval of the Cantor set. 
Since for any $t$ the set $\Pt(t)(E)$ has Lebesgue measure zero,  this implies that 
$$
\left|\bigcup_{t \in [0,1]\sm C}\Pt(t)(E\backslash\punct{f(t)})\right|=0,
$$
hence to prove \eqref{e:punchedwithf} it is enough to show
\begin{equation}\label{e:enough}
\left|\bigcup_{t \in C}\Pt(t)(E\backslash\punct{f(t)})\right|
<|P_0[E\sm U_0]|+\eps.
\end{equation}

By construction, for any $t\in C$ we have $U_{\phi_{n+1}(t)}\subset U_{\phi_{n}(t)}$.
By \eqref{e:fcenter} and the definition of $f$,
$$\bigcap_{n=1}^\infty U_{\phi_n(t)}\su \bigcap_{n=1}^\infty B_{\de_n}(\langle f_n(t)\rangle)\subset\langle f(t)\rangle.$$ 
Thus we have
\begin{align*}
    \bigcup_{t \in C} \Pt(t)(E\sm \punct{f(t)}) \su
    \bigcup_{n \in \N} \bigcup_{t \in C}\Pt(t)(E\sm {U_{\phi_n(t)}}).
\end{align*}
The right-hand side is an upward nested countable union of compact sets, hence it is Lebesgue measurable, and its Lebesgue measure is the limit of $|\cup_{t \in C}\Pt(t)(E\sm {U_{\phi_n(t)}})|$, so
for large enough $n$ we have
\begin{equation}\label{e:tCPtf}
\left| \bigcup_{t \in C} \Pt(t)(E\sm \punct{f(t)}) \right| \le
\left| \bigcup_{t \in C}\Pt(t)(E\sm {U_{\phi_n(t)}})\right| + \frac{\eps}2.
\end{equation}
Since $C\subset J^n$, this and \eqref{e:JnPtn} imply \eqref{e:enough}, which completes the proof.
%
%
\end{proof}

The following simple consequence might be interesting in itself.

\begin{corollary}\label{c:punchedsegment}
   Let $E$ be a closed line segment. 
Given some fixed $\varepsilon > 0$, $\eta>0$ and any (continuous) path $P_0:[0,1]\to\Iso$, there exists a path $\Pt:[0,1]\to\Iso$ with the same endpoints as $P_0$
within distance $\eta$ from a reparametrization of $P_0$,
and a continuous map $f:[0,1]\to E$ 
such that 
\begin{equation}\label{e:punchedwithfincor}
\left|\bigcup_{t \in [0,1]}\Pt(t)(E\sm\{f(t)\})\right| < \varepsilon.
\end{equation}
\end{corollary}

\begin{proof}
Any path is by definition continuous on $[a,b]$ and thus uniformly continuous. Thus any path can be approximated arbitrarily well with a polygonal path with the same endpoints, so we can suppose that $P_0$ is polygonal.
Then if we choose $U_0$ as the interior of $E$, we get $|P_0[E\sm U_0]|=0$,
hence Lemma~\ref{lemma:P exists} gives $\Pt$ with the required properties. 
\end{proof}

\begin{lemma}\label{the big guns}
    Let $E$ be a closed line segment (or circle), let $P_0:[0,1]\to\Iso$ be a polygonal path. Fix $\varepsilon >0$. Let $U_0$ be an open interval properly contained in the interior of 
$E$ (for the circle, two opposite open intervals). In the case of the segment, let $\tilde{E}$ be the extension of $E$ in both directions by the length of $U_0$, and in the case of the circle, simply let $\tilde{E} = E$. 
Let $c_0\in U_0$.
Then, there exists a path $P:[0,1]\to\Iso$ with the same endpoints as $P_0$ such that
    \begin{equation}\label{e:removedpoint}
        |P[E\setminus \punct{c_0}]| < |P_{0}[\tilde{E}\setminus U_{0}]|+\varepsilon.
    \end{equation}
    
Furthermore, in the linear case, if $U_0$ is in the middle third of $E$ then we can also guarantee
  \begin{equation}\label{e:onlyformiddlethird}  
    |P_{0}[E]\setminus P[\punct{c_0}]|< 2|P_0[\tilde{E}\setminus U_0] |+3\varepsilon.
    \end{equation}  
    
\end{lemma}
\begin{proof}
Applying Lemma~\ref{lemma:P exists} to $\Et$  we obtain
a path $\Pt:[0,1]\to\Iso$ with the same endpoints as $P_0$, 
and a continuous map $f:[0,1]\to \overline{U_0}$ (or simply $f:[0,1]\to E$ in the circular case)
such that
\begin{equation}\label{e:punchedEtwithf}
\left|\bigcup_{t \in [0,1]}\Pt(t)(\Et\backslash\punct{f(t)})\right| < |P_0[\Et\setminus U_0]| +\varepsilon
\end{equation}
and
\begin{equation}\label{e:P0EtminusPEt}
\left|P_0[\Et]\setminus \Pt[\Et]\right| <  \varepsilon.
\end{equation}

Extend $\Pt$ to $[-1,0]$ as $\Pt(0)$, to $[1,2]$ as $\Pt(1)$.
Let $f(-1)=f(2)=c_0$ and extend $f$ to $(-1,0)$ and $(1,2)$ 
by linear interpolation in the arc-length
coordinate. 
Then reparametrize the extended $\Pt$ and $f$ simultaneously to make the domain $[0,1]$. 
This way all the above properties of $\Pt$ and $f$ remain true and for the modified $f$ we also have $f(0)=f(1)=c_0$. 
Therefore we can suppose that $f(0)=f(1)=c_0$.

For any $t\in[0,1]$ let $T(t)\in\Iso$ be the translation (for the circle, rotation) preserving the line containing $E$ (or simply $E$) that maps $c_0$ to $f(t)$, and let $P(t)=\Pt(t)\circ T(t)$.
Then $P:[0,1]\to\Iso$ is continuous and it has endpoints $P_0(0)\circ T(0)$ and $P_0(1)\circ T(1)$.
Since $f(0)=f(1)=c_0$, $T(0)=T(1)=\id$, so $P$ has the same endpoints as $P_0$.
Note also that for any $t\in[0,1]$ we have
$$
P(t)(E\sm \punct{c_0})=\Pt(t)(T(t)(E\sm \punct{c_0}))=\Pt(t)(T(t)(E)\sm \punct{f(t)})\su \Pt(t)(\Et\sm \punct{f(t)}),
$$
hence by \eqref{e:punchedEtwithf} we obtain
$$
|P[E\setminus \punct{c_0}]| \le 
\left|\bigcup_{t \in [0,1]}\Pt(t)(\Et\backslash\punct{f(t)})\right|
< |P_{0}[\tilde{E}\setminus U_{0}]|+\varepsilon,
$$
which gives \eqref{e:removedpoint}. 

So it remains to prove \eqref{e:onlyformiddlethird} and we can suppose that $U_0$ is in the middle third of $E$. 
This condition implies that for any $t\in[0,1]$ we have $T(t)(E) \supset \overline{U_0}$, hence  $P[E]=(\Pt\circ T)[E]\supset \Pt[\overline{U_0}]$.
Thus 
$$\Pt[\tilde{E}]\setminus P[E] \su \Pt[\tilde{E}]\setminus \Pt[\overline{U_0}]\su \Pt[\tilde{E}\setminus \overline{U_0}],$$ 
so using also $E\su\Et$, \eqref{e:P0EtminusPEt}, $\punct{f(t)}\su \overline{U_0}$ and \eqref{e:punchedEtwithf} we obtain
\begin{align*}
    |P_0[E]\setminus  P[E]| &\leq |P_0[\tilde{E}]\setminus \Pt[\tilde{E}]| + |\Pt[\tilde{E}]\setminus  P[E]| \\
    & < {\varepsilon} +|\Pt[\tilde{E}\setminus \overline{U_0}]| \\
    &\leq \eps + \left|\bigcup_{t \in [0,1]}\Pt(t)(\tilde{E}\setminus \punct{f(t)})\right|  \\
    &< |P_0[\tilde{E}\setminus U_0]| + 2\varepsilon.
\end{align*}

Therefore, using this and the already proved \eqref{e:removedpoint} we get
\begin{align*}
    |P_0[E]\setminus  P[\punct{c_0}]| &\leq |P_0[E]\setminus P[E]| + | P[E]\setminus  P[\punct{c_0}]| \\
    &< \left(|P_0[\tilde{E}\setminus U_0]| + 2\varepsilon\right) + | P[E\setminus \punct{c_0}]| \\
    &< 2|P_0[\tilde{E}\setminus U_0]|+3\varepsilon
\end{align*}
which gives \eqref{e:onlyformiddlethird} and  completes the proof. 
\end{proof}

\bigskip

%

\section{Proofs of the main results}\label{s:main}

Now we are ready to prove the main results of the paper.

\begin{theorem}\label{t:bird}
    Let $E$ be a line segment, $q$ be any point of $E$,
$S$ be any bounded subset of the plane and $\eps>0$. 
Then 
there exists a path of orientation-preserving isometries of the plane, $P$, such that 
    \begin{equation}\label{tag1}
    \left|S \sm P[\{q\}]\right| < \varepsilon,
\end{equation} and 
\begin{equation}\label{tag2}
    \left|P[E\backslash\{q\}]\right| <\varepsilon.
\end{equation}

That is, the line segment $E$
can be moved continuously so that during the motion $E\sm\{q\}$ sweeps an area less than $\eps$, while the single point $q$
sweeps all points of $S$ but a set of area less than $\eps$.
\end{theorem}

\begin{proof}
Without loss of generality we can suppose that $E$ is closed and $q$ is the center of $E$ since otherwise we can extend $E$ to achieve this. 
Let $U_0$ be an open line segment inside the middle third of $E$ such that $q\in U_0$. 
Let $\tilde{E}$ be the line segment $E$ extended in both directions by the length of $U_{0}$, as in Lemma~\ref{the big guns}. 
We apply Cunningham's Bird Lemma (Theorem~\ref{Cunningham}) to get a polygonal path $P_{0}$ so that
\begin{equation}\label{e:P0U0coversS}
S\su P_{0}[U_{0}],
\end{equation}
 but 
\begin{equation}\label{e:restcoverssmall}
|P_{0}[\tilde{E}\setminus U_{0}]|<\varepsilon/{5}.
\end{equation} 

By applying Lemma~\ref{the big guns} 
we obtain a path $P$ such that 
\begin{equation}\label{e:Eminusq}
|P[E\sm\{q\}]|<|P_0[\tilde{E}\sm U_0]|+\eps/5
\end{equation}
and
\begin{equation}\label{e:P0EminusPq}
|P_0[E]\sm P[\{q\}]|<2|P_0[\tilde{E}\sm U_0]|+3\eps/5.
\end{equation}
Then \eqref{e:restcoverssmall} and \eqref{e:Eminusq} clearly imply \eqref{tag2}.
Since $U_0\su E$, \eqref{e:P0U0coversS} implies that $S\sm P[\{q\}]\subset P_0[E]\sm P[\{q\}]$, so \eqref{e:restcoverssmall} and \eqref{e:P0EminusPq} imply \eqref{tag1}, which completes the proof. 
\end{proof}

\begin{theorem}\label{t:strongKakeyaofcircleminustwopoints}
    The circle with two antipodal points removed has the strong Kakeya property. 
That is, for any $\rho_1,\rho_2\in\Iso$ and $\varepsilon>0$,
there is a path $P:[0,1]\to\Iso$ from $\rho_1$ to $\rho_2$ such that
\[
|P[E\sm\{c_0,c_1\}]|<\varepsilon,
\]
where $E$ is the circle and $c_0$ and $c_1$ are two antipodal points.    
\end{theorem}
\begin{proof}
Let $P_0: [0,1] \to\mathrm{Isom}^+(\R^2)$ be any segment of motion between $\rho_1$ and $\rho_2$. 
    Let $U_0 \su E$ be the union of two opposite open semicircles such that $c_0\in U_0$. Fix $\varepsilon>0$. Clearly, $|P_0[E\setminus U_0]| = 0$. Thus, by Lemma \ref{the big guns}, we attain a path $ P$ between $\rho_1$ and $\rho_2$ such that 
    $$| P[E\setminus \{c_0,c_1\}]| < |P_0[E \setminus U_0]| + \varepsilon = \varepsilon,$$
    as desired. 
\end{proof}

\begin{theorem}\label{t:punchedlineStrongKakeya}
A punched line in the plane has the strong Kakeya property with respect to any
absolutely continuous finite Borel measure.

More precisely, 
if $\mu$ is an absolutely continuous finite Borel measure on $\mathbb R^2$,
$\ell\su\R^2$ is a line and $p\in\ell$
then for any
$\rho_1,\rho_2\in\Iso$ and $\varepsilon>0$ there exists a path
$P:[0,1]\to\Iso$ from $\rho_1$ to $\rho_2$ such that
\[
\mu(P[\ell\sm\{p\}])<\varepsilon.
\]
%
\end{theorem}

\begin{proof}

Since $\mu$ is finite, for large enough $R$ we have $\mu(\R^2\sm D_R)<\eps/2$,
where $D_R$ denotes the open disc with radius $R$ and centered at the origin.
Let $E\su\ell$ be a long enough closed line segment such that $\rho_1(\ell\sm E)$ and $\rho_2(\ell\sm E)$ are disjoint from $D_{R+2}$.
Let $P_0:[0,1]\to\Iso$ be a path between $\rho_1$ and $\rho_2$ such that $P_0(t)(\ell\sm E)\su\R^2\sm D_{R+1}$ for every $t\in[0,1]$.

Since $\mu$ is absolutely continuous there exists $\de>0$ such that for any Borel set $B$ of Lebesgue measure less than $\de$ we have $\mu(B)<\eps/2$.
Applying Corollary~\ref{c:punchedsegment} with $\eps=\de$ and $\eta$ sufficiently small,
we get a path $\Pt:[0,1]\to\Iso$ between $\rho_1$ and $\rho_2$ and a continuous map $f:[0,1]\to E$ such that 
$\Pt(t)(\ell\sm E)\su \R^2\sm D_R$ for every $t\in[0,1]$ and 
\begin{equation*}
\left|\bigcup_{t \in [0,1]}\Pt(t)(E\sm\{f(t)\})\right| < \de.
\end{equation*}
Using also that $\mu(\R^2\sm D_R)<\eps/2$ we obtain
\begin{equation*}\label{e:punchedwithfincoragain}
\mu\left(\bigcup_{t \in [0,1]}\Pt(t)(\ell\sm\{f(t)\})\right) < \eps.
\end{equation*}

Let $T(t)$ be the translation parallel to $\ell$ that maps $p$ to $f(t)$ and for each $t\in[0,1]$ let $Q(t)=\Pt(t)\circ T(t)$.
Then $Q(t)(\ell\sm\{p\})=\Pt(t)(\ell\sm\{f(t)\})$, so we get $\mu(Q[\ell\sm\{p\}])<\eps$.

Now $Q(0)=\rho_1\circ T(0)$ and $Q(1)=\rho_2\circ T(1)$.
To join these maps with $\rho_1$ and $\rho_2$, respectively, let $L_1$ be a translational segment between $\rho_1$ and $\rho_1\circ T(0)$ and let
$L_2$ be a translational segment between $\rho_2\circ T(1)$ and $\rho_2$.
Then $L_i[\ell\sm\{p\}]\su \rho_i(\ell)$ ($i=1,2$), which have $\mu$-measure zero,
so by concatenating and reparametrizing $L_1$, $Q$ and $L_2$ we get a path $P$ with all the required properties.
%
\end{proof}


Now we use an idea from \cite{Cu}.  
Theorem~\ref{t:punchedlineStrongKakeya} can be easily made into an analogous theorem for the projective plane, and then using the standard two-to-one map from the sphere onto the projective plane we get the following corresponding result on a sphere for great circles with two antipodal points removed.
For completeness we give a proof without referring to projective geometry.

\begin{theorem}\label{t:greatcircle}
A great circle with two antipodal points removed has the strong Kakeya property on the sphere.

More precisely,
on a sphere $S\su\R^3$ if $W$ is a great circle with two antipodal points removed
then for any $\rho_1,\rho_2\in \IsoS$ and $\varepsilon>0$ there exists a
path $Q:[0,1]\to\IsoS$ such that
\[
Q(0)=\rho_1,\qquad Q(1)=\rho_2,
\qquad
\textrm{and}
\qquad
|Q[W]|<\varepsilon.
\]
%
\end{theorem}

\begin{proof}
We can clearly suppose that $S$ is the unit sphere centered at the origin and the removed points of $\rho_1(W)$ and $\rho_2(W)$ are not on the $z=0$ horizontal plane.
Let $U$ be the open upper hemisphere of $S$ and let $\Phi:U\to\R^2\times\{1\}$ be the projection from the origin. 
Then $\Phi$ is a homeomorphism and it also gives a homeomorphism between the set of open half great circles of $U$ and the set of lines of $\R^2\times\{1\}$.
Let $\mu$ be the pushforward of the surface measure on $U$ under the map $\Phi$.
Clearly, $\mu$ is a finite, absolutely continuous Borel measure on the plane $\R^2\times\{1\}$.

Fix an orientation of $W$, and let $p$ be the removed point of $W$
such that $\rho_1(p)\in U$. The sets $\Phi(\rho_1(W)\cap U)$ and $\Phi(\rho_2(W)\cap U)$
are oriented lines with one point removed. Hence, by
Theorem~\ref{t:punchedlineStrongKakeya}, there is a path
$P:[0,1]\to\Iso$ moving the first oriented punched line to the
second one such that the swept set has $\mu$-measure less than
$\varepsilon/2$.

Pulling this family back by $\Phi^{-1}$ and extending antipodally,
we obtain a continuous family $C_t$ of oriented great circles
with an antipodal pair of points removed. Mark the removed point
lying in $U$ by $p_t$.

An oriented great circle together with a marked point determines
uniquely an orientation-preserving isometry of $S$ relative to a
fixed such configuration, and this isometry depends continuously
on the data. Hence the family $(C_t,p_t)$ is induced by a continuous
path $\Qt:[0,1]\to\IsoS$ with $\Qt(0)=\rho_1$.

At the other endpoint we have either $\Qt(1)=\rho_2$ or
$\Qt(1)=\rho_2\circ J$, where $J$ is the half-turn about the axis
perpendicular to the plane of $W$. In the first case let $Q=\Qt$. 
In the latter case, we obtain $Q$ by concatenating
$\Qt$ with the corresponding rotation from $\rho_2\circ J$ to
$\rho_2$ and reparametrizing. 
During this additional motion, the punctured great circle remains in the same great circle, so the additional sweep has zero surface measure.

The part of the sweep in $U$ has surface measure less than
$\varepsilon/2$, its antipodal copy has the same measure, and the
points on the equator sweep a set of surface measure zero. Thus the
total swept surface measure is less than $\varepsilon$.
\end{proof}

If we combine Theorem~\ref{t:greatcircle} with a simple argument (which already appears in \cite{Cu}) we can give a full characterization of those subsets of a great circle that have the strong Kakeya property on a sphere.

\begin{corollary}
    Let $E$ be a great circle on a sphere $S$ centered at the origin. 
    Then $A\su E$ has the strong Kakeya property on the sphere if and only if $E\sm A$ contains a pair of antipodal points.
\end{corollary}

\begin{proof}
    The backward direction clearly follows from Theorem \ref{t:greatcircle}. 
    For the forward direction, let $\rho$ be a rotation of $\pi$ radians about an axis intersecting $E$ through the origin and suppose that $E\sm A$ contains no pair of antipodal points. 
    For any path $P\colon [0,1]\to \IsoS$ from the identity to $\rho$,  $P[E]$ is the entire surface of the sphere, for topological reasons. 
    So it suffices to show that the surface area of $P[A]$ is at least the surface area of $P[E\sm A]$.

If $p\in P[E\sm A],$ then $P(t)(x)=p$ for some $t\in[0,1], x\in E\sm A$. But since $E\sm A$ does not contain antipodes, $-x\in A$, so $-p=P(t)(-x)\in P[A].$ Thus $P[A]$ contains a reflected copy of $P[E\sm A],$ so we are done.
\end{proof}

\section{Applications: discontinuity of the area swept by continuously moving objects}\label{s:applications}

In this section we show how the results of Section~\ref{s:main} can be applied to obtain some discontinuity results about the area swept by continuously moving objects.

Let $E$ be a bounded subset of the plane and let $v\in\R^2$ be a given translation vector. For every path $P$ of orientation-preserving isometries, define $f_P:[0,1] \to \R$ by 
$$
f_P(t)=|P[E+vt]|,
$$
where $E+vt$ denotes the set $E$ translated by $vt$.

For fixed $E$ and $v$, we would like to decide whether $f_P$ is continuous for every path $P$.
It is easy to show that if $E$ is compact, then $f_P$ must be upper semicontinuous. 
    
The most natural case is when $E$ is just a point. By symmetry, we can assume $E$ to be the origin $(0,0)$. 
Perhaps unintuitively, Zs.~Fleiner \cite{Fl} constructed a path $P:[0,1]\to\Iso$ such that $|P[\{0\}]| = 1$, while $|P[\{x\}]| = 0$ for every point \mbox{$x\in [-1,1]^2 \setminus \{(0,0)\}$.} For this $P$, we clearly have that \[f_P(t)=|P[E+vt]|=|P[\{vt\}]|\] is not continuous at $t=0$ for every nonzero translation vector $v$. 
The construction is technically involved: it is an application of a highly nontrivial generalization of a construction due to Talagrand \cite{Ta}. In fact, we do not know of any simple direct proof of the discontinuity, even in this simple case.

We now show that Theorem~\ref{t:bird} easily gives another construction for any one point set $E=\{e\}$ and translation vector $v\in\R^2\sm\{(0,0)\}$.
Let $S$ be the unit square and $\varepsilon = \frac14$. Applying the theorem to the closed line segment $[e,e+v]$ with $q=e$, we obtain a path $P:[0,1]\to\Iso$ such that 
$|P[[e,e+v]\sm\{e\}]|<\frac14$ and $|P[\{e\}]| > \frac34$.
We then have
    $$f_P(0) = |P[\{e\}]| >\frac34,$$ and for all $t \in (0,1]$,
    $$f_P(t) \leq |P[[e,e+v]\setminus\{e\}]| < \frac{1}{4},$$
    showing that $f_P$ is discontinuous at $0$. 


In the remainder of this section we consider the cases when $E$ is either a closed or an open line segment.

    \begin{definition}
        A \textbf{rectangular map} is a function $[0,1] \to \R$ defined by
        $$ t\mapsto |P[E+vt]|,$$
        where $P$ is a path of orientation-preserving isometries, $E$ is a line segment (open or closed), and 
$v$ is a nonzero translation vector perpendicular to $E$. 
        The map is \emph{open} if $E$ is an open line segment and \emph{closed} if $E$ is a closed line segment.
    \end{definition}
\begin{definition}Given a line segment $E$, we say that a vector $v$ is \emph{diagonal} if it is neither parallel nor perpendicular to $E$.\end{definition}
    \begin{definition}
        A \textbf{diagonal map with translation vector $v$} is a function $[0,1] \to \R$ defined by 
                $$ t\mapsto |P[E+vt]|,$$
        where $P$ is a path of orientation-preserving isometries, $E$ is a line segment, and         
$v$ is a diagonal translation vector. The map is called open or closed based on whether $E$ is open or closed.
     \end{definition}

We are interested in whether all maps of a given type are continuous.
By invariance under similarities and reparametrization, this property
depends only on whether the segment is open or closed and on the
(unoriented) angle between the segment and the translation vector.
Note that Theorem~\ref{t:opendiagonal-intro} immediately follows from the following theorem and Theorem~\ref{t:bird}.

\begin{theorem}\label{t:opendiagonal}
    For any open line segment and any fixed diagonal vector $v$, there is an open diagonal map with translation vector $v$ that is discontinuous.
    
    In other words, if $E\su\R^2$ is an open line segment and $v\in\R^2$ is neither perpendicular nor parallel to $E$, then the plane can be moved continuously so that the area swept by $E+tv$ is not continuous in $t$. 
\end{theorem}
\begin{proof}
Let $\theta$ be the (unoriented) angle between $E$ and $v$.
Assume for contradiction that all open diagonal maps with angle
$\theta$ are continuous. Consider a line segment with endpoints
$A,C$, parallel to $E$, and midpoint $B$.
By Theorem \ref{t:bird}, let $P$ be a path of orientation-preserving isometries under which the closed segment from $A$ to $C$ with $B$ removed sweeps a set of measure less than $\varepsilon>0$ and $B$ sweeps a set of measure at least $1$.


Since open diagonal maps with 
angle $\theta$
are assumed to be continuous, for the open line segments $(A,B),(B,C)$ and sufficiently small $t$, we have that \[|P[(A,B)+tv]| <2\varepsilon, \]\[ |P[(B,C)+tv]| <2\varepsilon.\]Let $R$ be the reflection across the perpendicular bisector of $AC$, and let $v' = L_R(v)$, where $L_R$ denotes the linear part of $R$, as shown in Figure~\ref{fig:diagonal-translation}. 
Since $v'$ makes the same unoriented angle $\theta$ with the segment
as $v$, our continuity assumption applies to translation vector $v'$
as well.
Therefore, for sufficiently small $t$, we also have that \[|P[(A,B)+tv']| <2\varepsilon, \] \[ |P[(B,C)+tv']| <2\varepsilon.\]

Therefore, by choosing sufficiently small $t$, we obtain two new triples of points $A'= A +tv,B'= B +tv,C' = C+tv$ and $A'' = A + tv',B'' = B+ tv',C'' = C+ tv'$ such that under $P$ the open segments $(A',B')$, $(B',C')$, and $(A'',B'')$ sweep sets of measure less than $2\varepsilon$, while $(A',C')$ sweeps a set of measure at least $1-\varepsilon$. Also make sure the translation distance is small enough that the translated points appear in the order $A',A'',B',B'',C',C''$, as shown in Figure~\ref{fig:diagonal-translation}. Since $(A',C')$ sweeps a set of measure at least $1-\varepsilon$, but $(A',B')$ and $(B',C')$ each sweep a set of measure at most $2\varepsilon$, the single point $B'$ must sweep a set of measure at least $1-5\varepsilon$ under $P$. However, $B'$ lies between $A''$ and $B''$, so $(A'',B'')$ must also sweep a set of measure at least $1-5\varepsilon$, but we know $(A'',B'')$ sweeps a set of measure at most $2\varepsilon < 1-5\varepsilon$ (for $\varepsilon < \frac{1}{7}$), a contradiction.

\begin{figure}[ht]
\centering
\begin{tikzpicture}
\draw[blue, thick] (0,0) -- (0,3);
\draw[blue, thick] (0.6,-0.5) -- (0.6,3.5);

\draw[
    black,
    dashed,
    postaction={
        decorate,
        decoration={
            markings,
            mark=at position 0.75 with {\arrow{stealth}}
        }
    }
] (0,1.5) -- (0.6,2);
\draw[black,
    dashed,
    postaction={
        decorate,
        decoration={
            markings,
            mark=at position 0.75 with {\arrow{stealth}}
        }
    }] (0,1.5) -- (0.6,1);
\draw[black, dashed] (0,3) -- (0.6,3.5);
\draw[black, dashed] (0,3) -- (0.6,2.5);
\draw[black, dashed] (0,0) -- (0.6,0.5);
\draw[black, dashed] (0,0) -- (0.6,-0.5);

\node at (0.3,1.95) {$v$};
\node at (0.3,0.98) {$v'$};

\fill[blue] (0,0) circle (2pt) node[left] {$C$};
\fill[blue] (0,3) circle (2pt) node[left] {$A$};

\fill[blue] (0.6,-0.5) circle (2pt) node[right] {$C''$};
\fill[blue] (0.6,0.5) circle (2pt) node[right] {$C'$};
\fill[blue] (0.6,3.5) circle (2pt) node[right] {$A'$};
\fill[blue] (0.6,2.5) circle (2pt) node[right] {$A''$};

\fill[red] (0,1.5) circle (2pt) node[left] {$B$};
\fill[red] (0.6,2) circle (2pt) node[right] {$B'$};
\fill[red] (0.6,1) circle (2pt) node[right] {$B''$};
\end{tikzpicture}
\caption{$(A,C)$ translated by $tv$ and $tv'$ for small $t$.}
\label{fig:diagonal-translation}
\end{figure}
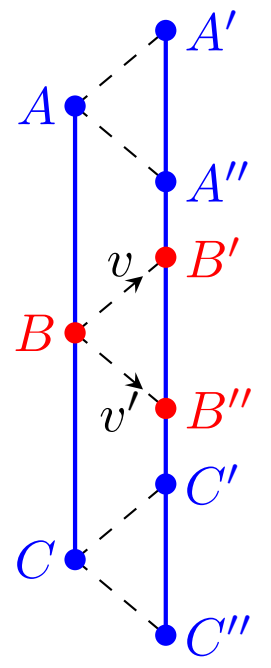
\end{proof}

It is now natural to think that discontinuity can occur in the other three cases (open rectangular, closed rectangular, and closed diagonal); however, 
we can prove only the following weaker result.

\begin{theorem}
If all open rectangular maps are continuous, then 
for any fixed closed segment and for every diagonal vector there is an associated closed diagonal map that is discontinuous.

In other words, if $E$ is an open line segment, $u$ is a nonzero vector perpendicular to $E$, and $v$ is a vector that is neither parallel nor perpendicular to $E$, then at least one of the following statements is false:
\begin{enumerate}
\item For any continuous motion of the plane, the area swept by $E+tu$ is continuous in $t$.
\item For any continuous motion of the plane, the area swept by $\overline{E}+tv$ is continuous in $t$.
\end{enumerate}
\end{theorem}

\begin{proof}
Assume for contradiction 
that all open rectangular maps are continuous, but there exist a closed line segment and a diagonal vector for which all associated closed diagonal maps are continuous. 
By invariance under similarities, it is clear that such a vector exists for \textit{any} closed line segment. Consider a line segment with endpoints $A,C$ and midpoint $B$. Let $v$ be a vector as above for the segment $[B,C]$. By Theorem~\ref{t:bird}, let $P$ be a path of orientation-preserving isometries such that, under $P$, the half-closed segment $[A,C)$ sweeps a set of measure less than $\varepsilon > 0$, and $C$ alone sweeps a set of measure at least $1$. For the remainder of the proof, diagonal translation is understood to mean translation in the direction of $v$. By the definition of $v$, translating $[B,C]$ diagonally by any sufficiently small distance, we get a new segment determined by points $B',C'$ such that under $P$, the segment $[B',C']$ sweeps a set of measure at least $|P[[B,C]]|-\varepsilon \geq 1-\varepsilon$. Also consider a perpendicularly translated copy of the segment $(A,C)$ defined by endpoints $A'',C''$ that lies on the same line as the diagonal translation. Assume the diagonal translation is small enough that the translated points appear in the order $A'',B',C',C''$, as shown in Figure~\ref{fig:segment-translations}. By continuity of open rectangular maps, if the perpendicular translation is small enough, then $(A'',C'')$ sweeps a set of measure at most $2\varepsilon$ under $P$, hence we can assume our initial diagonal translation is also small enough to guarantee this (since the lengths of the perpendicular and diagonal translations are proportional). Then
\begin{align*}
    1-\varepsilon &\leq |P[[B',C']]| \leq |P[(A'',C'')]| < 2\varepsilon,
\end{align*}
a contradiction for all sufficiently small $\varepsilon$.
\begin{figure}[ht]
\centering
\begin{tikzpicture}
  \draw[blue, thick] (0,0) -- (0,3);

  \draw[blue, thick] (0.45,0) -- (0.45,3);
  \draw[blue, thick] (0.45,0.7) -- (0.45,3);

  \draw[black, dashed] (0,1.5) -- (0.45,2.2);
  \draw[black, dashed] (0,0) -- (0.45,0.7);

  \draw[black, dashed] (0,3) -- (0.45,3);
  \draw[black, dashed] (0,0) -- (0.45,0);

  \fill[blue] (0,3) circle (2pt) node[left] {$A$};
  \fill[blue] (0,1.5) circle (2pt) node[left] {$B$};
  \fill[red]  (0,0) circle (2pt) node[left] {$C$};

  
  \fill[blue] (0.45,2.2) circle (2pt) node[right] {$B'$};
  \fill[red]  (0.45,0.7) circle (2pt) node[right] {$C'$};

  \fill[blue] (0.45,3) circle (2pt) node[right] {$A''$};
  \fill[red]  (0.45,0) circle (2pt) node[right] {$C''$};
\end{tikzpicture}
\caption{The diagonal translate of $[B,C]$ and the corresponding horizontal translate of $(A,C)$.}
\label{fig:segment-translations}
\end{figure}
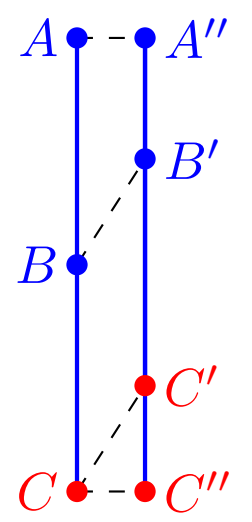
\end{proof}

As we remarked at the beginning of this section, it is easy to show that for compact $E$, the function $t\mapsto |P[E+tv]|$ is upper semicontinuous for any $v$, so any closed rectangular and closed diagonal maps are upper semicontinuous.
The following result shows that in the open case we have lower semicontinuity.

\begin{proposition}
    Open rectangular and open diagonal maps are lower semicontinuous.
\end{proposition}
\begin{proof}
We only prove the statement for open rectangular maps.
The proof for open diagonal maps 
can be done similarly.
 By symmetry, it is also enough to prove semicontinuity from one side and at one point; we will prove lower semicontinuity from the right at 0. Without loss of generality, we can assume that the segment is $((0,-1),(0,1))$.
    
    Since for any $P$,
    \[|P[((0,-1),(0,1))]| = \sup_{0<\eta<1} |P[[(0,-1+\eta),(0,1-\eta)]]|,\] it is enough to prove that for any fixed $\eta$ and $\varepsilon>0$ and all sufficiently small $\delta>0$ we have
    \begin{equation}\label{eq:lsc2}
    |P[I]\setminus P[((\delta,-1),(\delta,1))]|\le\varepsilon,
    \end{equation}
    where $I$ denotes the closed line segment $[(0,-1+\eta),(0,1-\eta)]$.

By dissecting the time interval into smaller intervals and proving the statement for each of them separately (with a sufficiently small $\varepsilon$), we can assume that $P$ is such that every point of $I$ remains within a distance of $\eta/2$ of its original position.
We claim that if $0<\delta_1<\delta_2$, then
\[ P[I] \cap P[((\delta_2,-1),(\delta_2,1))] \subset P[((\delta_1,-1),(\delta_1,1))].\] Let $L$ be the line $x = \delta_1$. If a point is contained in the left-hand side, then it must be swept by $L$ at some time $t_1$, and it is also swept by $I$ at some time $t_2$. Since each point of $L \sm ((\delta_1,-1),(\delta_1,1))$ is at a distance more than $\eta$ from $I$, a point swept by $L \sm ((\delta_1,-1),(\delta_1,1))$ cannot subsequently be swept by $I$. This proves the case $t_1 < t_2$; the other case is identical.
    The above shows that $P[I] \cap P[((\delta,-1),(\delta,1))]$ is an increasing family of sets as $\delta$ decreases.
    
    Let $S$ be the open half-strip $\{(x,y)\mid x>0, -1<y<1\}$. It is clear that \[\bigcup_{\delta>0} P[I] \cap P[((\delta,-1),(\delta,1))] = P[I] \cap P[S].\] Hence, by continuity of measure, if $\delta$ is small enough, then we have \begin{equation*}\label{eq:lsc3}
    |(P[I] \cap P[S]) \setminus P[((\delta,-1),(\delta,1))]| \le \varepsilon .\end{equation*}
    
Thus, to complete the proof of \eqref{eq:lsc2}, it is enough to show that  
   \begin{equation}\label{eq:lsc1}
        |P[I] \setminus P[S]| = 0.
    \end{equation}
Let $Q=P[I]\sm P[S]$ and take any $q\in Q$.
Then $q\in P(t_0)(I)$ for some $t_0$ and $P(t_0)(S)\su P[S]$, so $P(t_0)(S)\cap Q=\emptyset$.
Since the Lebesgue density of $q$ in $P(t_0)(S)$ is at least $\frac12$ , we obtain that the Lebesgue density of $q$ in $Q$ is at most $\frac12$.
Therefore, $Q$ cannot have a density point, hence the Lebesgue density theorem implies \eqref{eq:lsc1}.   
\end{proof}


\begin{thebibliography}{99}

\bibitem{Be}
A.~S.~Besicovitch, On Kakeya's problem and a similar one, \emph{Math. Z.} \textbf{27} (1928), 312--320.

\bibitem{CC}
A.~Chang and M.~Cs\"ornyei, 
The Kakeya needle problem and the existence of Besicovitch and Nikodym sets for rectifiable sets,
\emph{Proc. London Math Soc.} \textbf{118} (2019), 1084--1114.

\bibitem{CHL}
M.~Cs\"ornyei, K.~H\'era and M.~Laczkovich, 
Closed Sets with the Kakeya Property,
\emph{Mathematika} \textbf{63} (2017), 184--195.

\bibitem{Cu}
F.~Cunningham, Three Kakeya Problems, 
\emph{Amer.~Math.~Monthly} \textbf{81} (1974), 582--592.

\bibitem{Fe46}
H.~Federer,
Coincidence functions and their integrals,
\emph{Trans. Amer. Math. Soc.} \textbf{59} (1946), 441--466.

\bibitem{Fe52}
H.~Federer,
Measure and area,
\emph{Bull. Amer. Math. Soc.} \textbf{58} (1952), 306--378.

\bibitem{Fl}
Zs.~Fleiner, On Two Kakeya-type Problems, 
BSc thesis, 2025,
\url{https://www.math.elte.hu/thesisupload/thesisfiles/2025bsc_mat3y-am57lq.pdf}

\bibitem{HL}
K.~H\'era and M.~Laczkovich, The Kakeya Problem for Circular Arcs,
\emph{Acta Math. Hungar.} \textbf{150} (2016), 479--511.

\bibitem{Ta}
M. Talagrand, Sur la mesure de la projection d’un compact et certaines familles
de cercles, 
\emph{Bulletin des Sciences Mathématiques} \textbf{104} (1980), 225--231.
\end{thebibliography}
\end{document}